\documentclass[pdftex]{amsart}
\usepackage{subfigure}
\usepackage{euscript}
\usepackage[all]{xy} 
\usepackage{varioref}
\usepackage{lscape}
\usepackage{ifthen}
\usepackage{ifpdf}
\ifpdf
\newcommand{\drawfile}[1]%
{\includegraphics{{#1}.pdf}}
\else
\usepackage{epsfig}
\newcommand{\drawfile}[1]%
{\epsfig{file={#1}.eps}}
\fi
\usepackage{multirow}
\newtheorem{theorem}{Theorem} 
\newtheorem{corollary}{Corollary}
\newtheorem{lemma}{Lemma} 
\newtheorem{proposition}{Proposition}

\theoremstyle{remark}
\newtheorem{definition}{Definition} 
\newtheorem{example}{Example}
\newtheorem{remark}{Remark}

\newcommand{\LieDer}{\ensuremath{\EuScript L}}
\newcommand{\hook}{\ensuremath{\mathbin{ \hbox{\vrule height1.4pt
        width4pt depth-1pt \vrule height4pt width0.4pt depth-1pt}}}}

\newcommand{\R}[1]{\ensuremath{\mathbb{R}^{#1}}}
\newcommand{\C}[1]{\ensuremath{\mathbb{C}^{#1}}}

\newcommand{\CP}[1]{\ensuremath{\mathbb{CP}^{#1}}}
\newcommand{\RP}[1]{\ensuremath{\mathbb{RP}^{#1}}}

\newcommand{\SpinorPlus}[1]{
\ensuremath{
\mathbb{S}^{+}\left({#1}\right)
}
}
\newcommand{\SpinorMinus}[1]{
\ensuremath{
\mathbb{S}^{-}\left({#1}\right)
}
}
\newcommand{\SO}[1]{\operatorname{SO}\left({#1}\right)}

\newcommand{\PO}[1]{\ensuremath{\mathbb{P}O\left({#1}\right)}}

\newcommand{\GL}[1]{\operatorname{GL}\left({#1}\right)}

\newcommand{\gl}[1]{\mathfrak{gl}\left({#1}\right)}
\newcommand{\SL}[1]{\operatorname{SL}\left({#1}\right)}
\newcommand{\PSL}[1]{\mathbb{P}\SL{#1}}

\newcommand{\slLie}[1]{\mathfrak{sl}\left({#1}\right)}

\newcommand{\Symp}[1]{\ensuremath{\operatorname{Sp}\left({#1}\right)}}

\newcommand{\Sym}[2]{\ensuremath{\operatorname{Sym}^{#1}\left({#2}\right)}}
\newcommand{\Lm}[2]{\ensuremath{\Lambda^{#1} \left ( {#2} \right )}}
\newcommand{\nForms}[2]{\ensuremath{\Omega^{#1} \left ( {#2} \right
    )}}   \DeclareMathOperator{\Ad}{Ad}

\newcommand{\trans}[1]{{}^t{#1}}

\newcommand{\Lag}[1]{\ensuremath{\operatorname{Lag}\left(#1\right)}}
\newcommand{\Norm}[1]{\nu{#1}}
\newcommand{\Conorm}[1]{\nu^*{#1}}
\newcommand{\Gtot}{\ensuremath{G}} 
\newcommand{\gtot}{\mathfrak{g}}
\newcommand{\Xpt}[1]{\ensuremath{{#1}_0}}
\newcommand{\Xgeod}[1]{\ensuremath{{#1}_1}}
\newcommand{\Xptgeod}[1]{\ensuremath{{#1}_+}}
\newcommand{\Gpt}{\Xpt{\Gtot}} 
\newcommand{\gpt}{\Xpt{\gtot}}
\newcommand{\Ggeod}{\Xgeod{\Gtot}} 

\newcommand{\Gptgeod}{\Xptgeod{\Gtot}}

\newcommand{\parab}{\ensuremath{\mathfrak{p}}}
\newcommand{\Parab}{\ensuremath{P}}

\newcommand{\Proj}[1]{\mathbb{P}^{#1}}
\newcommand{\OO}[1]{
  \ensuremath{
    \mathcal{O}
    \ifthenelse{\equal{#1}{0}}
      {}
      {\left({#1}\right)}
  }
}
\newcommand{\OOp}[2]{
  \ensuremath{
    \mathcal{O}
    \ifthenelse{\equal{#1}{0}}
      {}
      {\left({#1}\right)}
    \ifthenelse{\equal{#2}{1}}
      {}
      {^{\oplus{#2}}}
  }
}

\newcommand{\semiprol}[2]
{
\ensuremath{
{#1}^{\left<{#2}\right>}
}
}
\newcommand{\prol}[2]
{
\ensuremath{
{#1}^{\left({#2}\right)}
}
}
\newcommand{\om}
{
\ensuremath{\omega}
}
\newcommand{\omsb}
{
\ensuremath{\bar{\omega}}
}

\newcommand{\s}[1]{\slLie{2,\C{}}_{#1}}
\newcommand{\Sl}[1]{\SL{2,\C{}}_{#1}}

\newcommand{\quot}[2]{#1/#2}

\newcommand{\al}[1]{#1}
\newcommand{\mn}[1]{}

\begin{document}
\title{Complete complex parabolic geometries}
\author{Benjamin McKay}
\address{School of Mathematical Sciences \\
University College Cork \\
  Cork, Ireland 
} \email{B.McKay@UCC.ie}
\date{\today} 
\thanks{Thanks to Robert Bryant, Robert Gunning and Claude LeBrun, and
especially to J. M. Landsberg for his insight
which made vital contributions to this paper.}
\begin{abstract}
Complete complex parabolic geometries
(including projective connections
and conformal connections) are flat
and homogeneous. This is the first global theorem
on parabolic geometries.
\end{abstract}
\maketitle
\tableofcontents
\section{Introduction}
Cartan geometries have played a fundamental role 
in unifying differential geometry (see Sharpe \cite{Sharpe:1997})
and in the theory of complex manifolds, particularly
complex surfaces (see Gunning \cite{Gunning:1978}).\mn{Removed ``I will prove''.}
\begin{theorem}
The only complete complex
Cartan geometry modelled on a rational homogeneous
space is the standard flat Cartan geometry on the
ratonal homogeneous space. In particular 
the two most intensely studied types of complex Cartan
geometries, projective connections and conformal
connections, when complete are flat. 
\end{theorem}
\par\noindent{}This is the first global theorem
on parabolic geometries. It is also the first characterization
of rational homogeneous varieties among complex manifolds 
which doesn't rely on symmetries. But much more importantly, this
theorem is the crucial ingredient to solve
the global twistor problem for complex parabolic geometries:
when can one apply a twistor transform to a parabolic geometry
(see McKay \cite{McKay:2006} for that problem and its solution). The theorem is very surprising
(as attested to by numerous geometers), especially by comparison
with Riemannian geometry where completeness is an automatic consequence of
compactness. The proof uses only classical techniques invented by Cartan.

The concept of completeness of Cartan geometries
is subtle, raised explicitly for the first time by
Ehresmann \cite{Ehresmann:1936} (also see 
Ehresmann \cite{Ehresmann:1938,Ehresmann:1951},
Kobayashi \cite{Kobayashi:1954}, 
Kobayashi \& Nagano \cite{KobayashiNagano:1964},
Clifton \cite{Clifton:1966}, 
Blumenthal \cite{Blumenthal:1987}, 
Bates \cite{Bates:2004}),
and plays a central role in the book of Sharpe \cite{Sharpe:1997},
but is also clearly visible beneath the surface in numerous works of
Cartan. Roughly speaking, completeness concerns
the ability to compare a geometry to some notion of
flat geometry, by rolling along curves. It is a
transcendental concept, loosely analoguous to
positivity conditions which are more familiar to algebraic geometers.
Indeed the results proven here are very similar
in spirit to those of Hwang \& Mok \cite{Hwang/Mok:1997}
and Jahnke \& Radloff \cite{Jahnke/Radloff:2002,Jahnke/Radloff:2004,Jahnke/Radloff:2004b},
but with no overlap; those authors  
work on smooth projective algebraic varieties,
and employ deep Mori theory to produce rational
circles (we will define the term \emph{circle} shortly), while
I take advantage of completeness to build rational
circles, working on any (not necessarily closed or K\"ahler)
complex manifold, more along the lines
of Hitchin \cite{Hitchin:1982}. They show that
minimal rational curves are circles, while I show
that complete circles are rational. We both start
off with a Cartan connection, and show that it is flat;
however, my Cartan connections include
all parabolic geometries, while \al{their Cartan connections} include
only those modelled on cominiscule varieties.

\section{Cartan connections}
We will need some theorems about Cartan geometries, mostly drawn from
Sharpe \cite{Sharpe:1997}.
\begin{definition}
  If $\Gpt \subset \Gtot$ is a closed subgroup of a Lie group and $\Gamma
  \subset \Gtot$ is a discrete subgroup, call the (possibly singular)
  space $\Gamma \backslash \Gtot/\Gpt$ a \emph{double coset space}.  If
  $\Gamma \subset \Gtot$ acts on the left on $\Gtot/\Gpt$ freely and properly,
  call the smooth manifold $\Gamma \backslash \Gtot/\Gpt$ a \emph{locally Klein
    geometry}.
\end{definition}
\begin{definition}
  A \emph{Cartan pseudogeometry} on a manifold $M$, modelled on a
  homogeneous space $\Gtot/\Gpt$, is a principal right $\Gpt$-bundle $E \to M$,
  (with right $\Gpt$ action written $r_g : E \to E$ for $g \in \Gpt$),
  with a 1-form $\om \in \nForms{1}{E} \otimes \gtot$,
  called the \emph{Cartan pseudoconnection} (where $\gtot, \gpt$ are
  the Lie algebras of $\Gtot,\Gpt$), so that $\om$ identifies each tangent
  space of $E$ with $\gtot$.  For each $A \in \gtot$, 
  let $\vec{A}$ be the vector field on $E$ satisfying $\vec{A} \hook
  \om = A$.  A Cartan pseudogeometry is called a \emph{Cartan
    geometry} (and its Cartan pseudoconnection called a \emph{Cartan
    connection}) if  
    \[
    \vec{A} = \left. \frac{d}{dt}r_{e^{tA}}\right|_{t=0}
    \]
    for all $A \in \gpt$, and $r_g^* \om = \Ad_g^{-1} \om$ for all $g \in \Gpt$.
\end{definition}
\begin{proposition}[Sharpe \cite{Sharpe:1997}]
  A locally Klein geometry is a Cartan geometry with the left
  invariant Maurer--Cartan 1-form of $\Gtot$ as Cartan connection.
\end{proposition}
\begin{definition}
  If $M_0 \to M_1$ is a local diffeomorphism, any Cartan geometry on
  $M_1$ pulls back to one on $M_0$, by pulling back the bundle and
  Cartan connection.
\end{definition}
\begin{definition}
  If $\Gamma \backslash \Gtot/\Gpt$ is a double coset space, and $\Gamma$ has
  a discontinuous action on a manifold $M$, commuting with a local
  diffeomorphism $M \to \Gtot/\Gpt$, then define a Cartan geometry on
  $\Gamma \backslash M$ given by pulling back via $M \to \Gtot/\Gpt$ and then
  quotienting by the automorphisms of $\Gamma$. \al{This is called a
  \emph{quotient of a pullback}.}
\end{definition}
\begin{definition}
A Cartan geometry modelled on $\Gtot/\Gpt$ is \emph{flat} if it is locally
isomorphic to $\Gtot/\Gpt$.
\end{definition}
\begin{theorem}[Sharpe \cite{Sharpe:1997}]
  Every flat Cartan geometry on a connected manifold is the quotient
  $M = \Gamma \backslash \tilde{M}$ of a pullback $\tilde{M} \to \Gtot/\Gpt$,
	with $\Gamma \subset \Gtot$ \al{being} a subgroup of $\Gtot$, acting as deck transformations on $\tilde{M}$.
\end{theorem}
\begin{proof}
Let $E \to M$ be any flat Cartan geometry, with Cartan
connection $\om$.
Put the differential system $\om=g^{-1} \, dg$
on the manifold $E \times \Gtot$.
By the Frobenius theorem, $E \times \Gtot$
is foliated by leaves (maximal connected
integral manifolds). 

We have the diagonal right action of $\Gpt$
on $E \times \Gtot$, and the left
action $g' \left(e,g\right)=\left(e,g'g\right)$.
Because the system
is invariant under both actions, these
actions permute leaves. 

Define vector fields $\vec{A}$
on $E \times \Gtot$
by adding the one from 
$E$ with the one (by the same name)
from $\Gtot$. The flow of $\vec{A}$ on $\Gtot$
is defined for all \al{times}, so
the vector field $\vec{A}$
on $E \times \Gtot$
has flow through a point 
$(e,g)$ defined for as long 
as the flow is defined down on $E$.
These vector fields $\vec{A}$ for $A \in \gtot$ 
satisfy $\om=g^{-1} \, dg$, so the leaves 
contain the flow lines of the $\vec{A}$.

If $L$ is a leaf, then the set $\Lambda = G_0 L$
(under diagonal $\Gpt$-action)
is a union of disjoint
leaves and thereby an immersed submanifold. 
Call such a set $\Lambda = G_0 L$ a
\emph{key}.
Let us show now that $\Gtot$ acts
transitively on keys.
Let $\Lambda_0$ and $\Lambda_1$ be
two keys containing points
$\left(e_i,g_i\right) \in \Lambda_i$.
Since $M$ is connected, after \al{$\Gpt$-action
arranges that $e_0$ and $e_1$ belong
to the same path component of $E$, we can} draw a 
path
from $e_0$ to $e_1$ in $E$, consisting
of finitely many flows of $\vec{A}$ vector 
fields; clearly such a path lifts to
a path in $\Lambda_0$ from $\left(e_0,g_0\right)$ 
to
some point $\left(e_1,g_1'\right)$ so after 
$\Gtot$-action
we get to $\left(e_1,g_1\right) \in \Lambda_1$. 
Therefore $\Gtot$ acts transitively on keys.

We can define the holonomy group $H_{e_0}$
of each point $e_0 \in E_0$ as the set
of all $g \in G$ so that $g \Lambda = \Lambda$
for $\Lambda$ the key through $\left(e_0,1\right)$.

Pick a key $\Lambda$.
The inclusion
$\Lambda \subset E \times \Gtot$
defines two local diffeomorphisms
$\Lambda \to E$ and $\Lambda \to \Gtot$.
\al{Both are} $\vec{A}$-equivariant. Consider
the first of these \al{diffeomorphisms}. Let $F$ be 
a fiber
of $\Lambda \to E$ over some point $e \in E$.
Define local coordinates on $E$ by
inverting the map
$A \in \gtot \mapsto e^{\vec{A}} e \in E$
near $A=0$. (This map is only defined
near $A=0$, and is a diffeomorphism
in some neighborhood, say $U$, of $0$). Then map 
\[
U \times F \to \Lambda 
\] 
by $(A,f) \mapsto e^{\vec{A}}f$, clearly a local
diffeomorphism. Therefore $\Lambda \to E$
is a covering map, and $\Gpt$-equivariant,
so covers a covering map $\tilde{M}=\Lambda/\Gpt 
\to M=E/\Gpt$.
Thus $\Lambda \to \tilde{M}$ is the pullback
bundle of $E \to M$. 

From now on we 
will pick a point $m_0 \in M$
and a point $e_0 \in E_{m_0}$, and
let $\tilde{E}$ be the key passing 
through the point $\tilde{e}_0=\left(e_0,1\right) \in E \times G$.

The map $\tilde{E} \to \Gtot$ is 
$\Gpt$-equivariant, 
so descends to a map $\tilde{M} \to \Gtot/\Gpt$.
On $\tilde{E}$, $\om=g^{-1} \, dg$,
so this map is a pullback of Cartan geometries.

For each absolutely continuous
loop $m(t)$ in $M$,
$0 \le t \le 1$,
say with $m(0)=m(1)=m_0$,
lift to a path $e(t) \in E$ and a path
with $e(0)=e_0$. There is a unique
lift of $e(t)$ to a path
\(
\tilde{e}(t) \in \tilde{E}
\)
with $\tilde{e}(0)=\tilde{e}_0$.
The path $e(t)$ must end at some point
$e(1)=e_0 h$, for a unique $h \in \Gpt$.
We can then write
\[
\tilde{e}(1)=\left(e_0h,gh\right),
\]
for a unique $h \in \Gpt$ and $g \in \Gtot$.
Let $\gamma$ be the homotopy class of
the loop $m(t)$ in $\pi_1(M)$.
We define
\[
\operatorname{hol}\left(\gamma\right)=g.
\]

If we have a piecewise continuously
differentiable homotopy of paths,
say $m(s,t)$ (thought of as a path for each
fixed value of $s$) and lift
to an arbitrary family of paths $e(s,t)$,
say from $e(s,0)=e_0$ to some point $e_0 h(s)$,
for each $s$, then 
the lift to $\tilde{E}$ is a family
of paths $\tilde{e}(s,t)$,
say with $\tilde{e}(s,0)=\tilde{e}_0$.
The other end point is somewhere
\[
\tilde{e}(s,1)=\left(e_0h(s),g(s)h(s)\right).
\]
Lying in $\tilde{E}$ ensures that this
path satisfies
\[
\frac{d}{ds}e_0 h(s) \hook \omega = 
\left(g(s)h(s)\right)^{-1} 
\frac{d}{ds}\left(g(s) h(s)\right),
\]
which expands out to say
\[
h^{-1} \frac{dh}{ds} = h^{-1} \frac{dh}{ds}
+ g^{-1}\frac{dg}{ds},
\]
or in other words $g(s)$ is constant.
Therefore $\operatorname{hol}(\gamma)$
is well defined. It is then easy
to prove that 
\[
\operatorname{hol} \colon
\pi_1(M) \to H_{e_0}
\]
is a group morphism.
Let $\Delta : \tilde{E} \to G$ be
the map $\Delta(e,g)=g$. 
Let $\pi_1(M)$ act on $\tilde{E}$
by
\[
\gamma(e,g)=\left(e,
\operatorname{hol}(\gamma)g\right).
\]
The map $\Delta$ is equivariant
under the $\pi_1(M)$-action,
the right $\Gpt$-action, and the
left $\Gtot$-action.
The $\pi_1(M)$-action commutes with the right
$\Gpt$-action, so descends to
an action on $\tilde{M}$, transitive
on the fibers of $\tilde{M} \to M$,
giving the deck transformations.
We can define the \emph{developing
map}
\[
\delta \colon 
\left(\tilde{M},\tilde{m}_0\right)
\to
\left(\Gtot/\Gpt,1 \cdot \Gpt\right),
\]
as the quotient of $\Delta$ by the
right $\Gpt$-action. Then clearly
\[
\delta \circ \gamma = 
\operatorname{hol}\left(\gamma\right)
\circ \delta
\]
for all $\gamma \in \pi_1(M)$.
\end{proof}
Warning: by definition, $\Gamma \backslash \tilde{M}$ will
be smooth, but $\Gamma \backslash \Gtot/\Gpt$ might not be.
\begin{definition}[Ehresmann \cite{Ehresmann:1936}]
  A Cartan geometry is \emph{complete} if all of the $\vec{A}$ vector
  fields are complete.
\end{definition}
For example, if the Cartan geometry is a pseudo-Riemannian
geometry, then this definition of completeness is equivalent to the usual
one (see Sharpe \cite{Sharpe:1997}).
\begin{lemma}\label{lemma:InfSymm}
If a Cartan geometry is complete, then all infinitesimal
symmetries of the Cartan geometry are complete vector
fields.
\end{lemma}
\begin{proof}
Infinitesimal symmetries are vector fields $X$ on $E$
for which $r_{g*} X = X$ and $\LieDer_X \om=0$.
The vector fields $\vec{A}$ thus satisfy $\left[\vec{A},X\right]=0$.
Therefore the flows lines of $X$ are permuted by
the flows of $\vec{A}$, so that the time for which
the flow is defined is locally constant. But then
that time can't diminish as we move along the flow, 
so $X$ must be complete.
\end{proof}
\begin{lemma}
  Pullback via covering maps preserves and reflects completeness.
\end{lemma}
\begin{corollary}[Sharpe \cite{Sharpe:1997}]\label{cor:CompleteGivesKlein}
  A Cartan geometry on a connected manifold is complete and flat just
  when it is locally Klein.
\end{corollary}
Warning: If the Cartan geometry we start with is modelled on $\Gtot/\Gpt$,
then as a locally Klein geometry it might be modelled on $\Gtot'/\Gpt$, for
some Lie group $\Gtot'$ with the same Lie algebra $\gtot$ as
$\Gtot$, and which also contains $\Gpt$ as a closed subgroup,
with the same Lie algebra inclusion $\gpt \subset
\gtot$; we can assume that $\Gtot'/\Gpt$ is connected.
\begin{corollary}[Sharpe \cite{Sharpe:1997},
  Blumenthal \& Hebda \cite{BlumenthalHebda:1989}]\label{cor:mutual}
  If two complete flat Cartan geometries on connected manifolds have
  the same model $\Gtot/\Gpt$, then there is a complete flat Cartan geometry
  on a mutual covering space, pulling back both.  In particular, if
  both manifolds are simply connected, then their Cartan geometries are isomorphic.
\end{corollary}
\begin{proof}
  Without loss of generality, $M_j = \Gtot_j/\Gpt$, \al{with} $\Gtot_j$ 
  playing the role of
  $\Gtot'$ in the previous corollary. Build a Lie
  group $\tilde{\Gtot}$ \al{covering both} $\Gtot_j$, and containing $\Gpt$.
\end{proof}
\begin{definition}
  A group $\Gamma$ \emph{defies} a group $\Gtot$ if every morphism $\Gamma
  \to \Gtot$ has finite image.  For example, if $\Gamma$ is finite, or $\Gtot$
  is finite, then $\Gamma$ defies $\Gtot$.
\end{definition}
I have only one theorem to add to Sharpe's collection:
\begin{theorem}\label{thm:compactToFlat}
  A flat Cartan geometry, modelled on $\Gtot/\Gpt$, defined 
  on a compact connected base
  manifold $M$ with fundamental group defying $\Gtot$, is a locally Klein
  geometry.
\end{theorem}
\begin{proof}
  Write $M = \Gamma \backslash \tilde{M}$ for some pullback $\tilde{M}
  \to \Gtot/\Gpt$.  $\Gamma = \pi_1(M)/\pi_1(\tilde{M}) \subset
  \Gtot$, and since  $\pi_1(M)$ defies $\Gtot$, $\Gamma$ is finite, 
  and $\tilde{M}$ is
  compact. The local diffeomorphism $\tilde{M} \to \Gtot/\Gpt$ must therefore
  be a covering map to its image. The bundles on which the
  Cartan connections live, say
  \[
  \xymatrix{
    E \ar[d] & \tilde{E} \ar[l] \ar[d] \ar[r] & \Gtot \ar[d] \\
    M        & \tilde{M} \ar[l] \ar[r] & \Gtot/\Gpt \\
  }
  \]
  are all pullbacks via covering maps, so completeness is preserved
  from $\Gtot/\Gpt$ and reflected to $M$.
\end{proof}
\al{Let us see that this theorem brings new insight.}
\begin{corollary}
  If $\dim \Gtot/\Gpt \ge 4$, then 
  infinitely many compact manifolds bear no flat
  $\Gtot/\Gpt$-Cartan geometry.
\end{corollary}
\begin{proof}
  Construct manifolds with fundamental group defying $\Gtot$, following
  Massey \cite{Massey:1967}. The fundamental group can be finite or
  infinite, as long as it has no quotient group belonging to $\Gtot$. For
  example, the fundamental group could be a free product of finitely presented
  simple groups not belonging to $\Gtot$.
\end{proof}
\subsection{Mutation}
\begin{definition}
  A locally Klein geometry $\Gamma \backslash \Gtot'/\Gpt$,
  where $\Gtot'$ is a Lie group containing $\Gpt$ as a closed subgroup 
  and of the same dimension
  as $\Gtot$ (but perhaps with a different Lie algebra)
  is called a \emph{mutation} of $\Gtot/\Gpt$.
\end{definition}
The standard example of a mutation is to take $\Gtot$ the
group of rigid motions of the Euclidean plane, $\Gpt$ the
stabilizer of a point, and let $\Gtot'$ be the group
of rigid motions of a sphere or of the hyperbolic plane.
The sphere and hyperbolic plane are then viewed as 
``Euclidean planes of constant curvature,'' 
although one could just as well view the Euclidean plane
as a ``sphere of constant negative curvature''; see
Sharpe \cite{Sharpe:1997} for more on this picture.
\begin{lemma}
Any mutation of a homogeneous space is complete.
\end{lemma}
\begin{proof}
The vector fields $\vec{A}$ are precisely
the left invariant vector fields on $\Gtot'$,
so they are complete.
\end{proof}
\begin{definition}
  Say that a connected homogeneous space $\Gtot/\Gpt$ is \emph{immutable} 
  if every connected mutation of it is locally isomorphic (and hence a locally Klein geometry 
  $\Gamma \backslash \tilde{\Gtot}/\Gpt$ for some covering Lie group
  $\tilde{\Gtot} \to \Gtot$).
\end{definition}
By contrast to Euclidean geometry, we will see 
(theorem~\vref{thm:Immutable}) that real and complex 
rational homogeneous spaces are immutable.
\subsection{Curvature}
\begin{lemma}[Sharpe \cite{Sharpe:1997} p.188]
  Take a Cartan geometry $E \to M$ with Cartan connection $\om$,
  modelled on $\Gtot/\Gpt$. Write the curvature as
\begin{align*}
  \nabla \om &= d \om + \frac{1}{2}
  \left[\om,\om\right] \\
  &= \frac{1}{2} \kappa \omsb \wedge \omsb,
\end{align*}
where $\omsb = \om + \gpt \in \nForms{1}{\om}\otimes
\left( \gtot/\gpt \right)$ is the \emph{soldering form} and $\kappa :
E \to \gtot \otimes
\Lm{2}{\gtot/\gpt}^*$.  Then under right
$\Gpt$ action $r_g : E \to E$,
\[
r_g^* \om = \Ad_g^{-1} \om, r_g^* \nabla \om = \Ad_g^{-1} \nabla \om.
\]
\end{lemma}
\begin{lemma}
A Cartan geometry is flat just when $\kappa=0$.
\end{lemma}
\begin{proof}
Apply the Frobenius theorem to the equation $\om = g^{-1} \, dg$
on $E \times \Gtot$; the resulting foliation of $E \times \Gtot$
consists of the graphs of local isomorphisms.
\end{proof}
\begin{theorem}[Sharpe \cite{Sharpe:1997}]\label{thm:ConstantCurv}
  A Cartan geometry is complete with
  constant curvature $\kappa$ just when it is a mutation.
\end{theorem}
\begin{proof}
	Taking the exterior derivative of the structure
	\al{equation} 
    \[d \om = - \frac{1}{2}\left[\om,\om\right]+\frac{1}{2} \kappa \omsb \wedge \omsb, 
    \]
    shows that the Lie bracket
    \[
    \left[A,B\right]'=\left[A,B\right]+\kappa(A,B)
    \]
	satisfies the Jacobi identity, so is a Lie bracket
	for a Lie algebra $\mathfrak{g}'$, containing
	$\mathfrak{g}_0$ as a subalgebra. The trouble is
	that perhaps every Lie group $\Gtot'$ with that Lie
	algebra fails to contains $\Gpt$.
		
	Without loss of generality, assume $\pi_1\left(M\right)=1$.
	The vector fields $\vec{A}$ satisfy 
    $\left[\vec{A},\vec{B}\right]=\overrightarrow{[A,B]'}$,
    so by Palais' theorem \cite{Palais:1957}
    generate a Lie group action, of a Lie group
    locally diffeomorphic to $E$. Let $\Gtot'$
    be the group of Cartan geometry automorphisms of $E$.
    Because the Cartan geometry is complete, the
    infinitesimal symmetries are complete by \al{Lemma}~\vref{lemma:InfSymm}.
	So the Lie algebra of $\Gtot'$ is precisely the
	set of infinitesimal symmetries. Consider
	the sheaf of local infinitesimal symmetries.
	The monodromy of this sheaf as we move around
	a loop in $E$, looking at the exact sequence
	\[
	\pi_1\left(\Gpt\right) \to \pi_1\left(E\right) \to \pi_1\left(M\right) = 1,
	\]
	comes from the monodromy around loops in $\pi_1\left(\Gpt\right)$.
	But local infinitesimal symmetries are invariant under
	flows of $\vec{A}$, so are invariant under the
	identity component of $\Gpt$, hence under $\pi_1\left(\Gpt\right)$.
	Therefore the local infinitesimal symmetries
	(coming from the local identification of $E$
	with a Lie group) extend globally, ensuring
	that $\Gtot'$ is locally diffeomorphic to $E$
	and acts locally transitively on $E$.
	Moreover, the vector fields $\vec{A}$
	pull back to the left invariant
	vector fields on $\Gtot'$.
	
	The infinitesimal symmetries induce
	vector fields on $M$, because they
	are $\Gpt$-invariant, and they are
	locally transitive on $M$. Since
	$M$ is connected, $\Gtot'$ acts
	transitively on $M$. 
	If $g' \in \Gtot'$ fixes some point $e_0 \in E$
	then commuting with $\vec{A}$ flows ensures that it fixes all nearby
	points of $E$, and therefore fixes all nearby
	points of $M$, and so fixes all points of $M$,
    so fixes all points of $E$. Therefore $\Gtot'$
    acts on $E$ locally transitively and
    without fixed points.

	Write $p_j : \left(e_0,e_1\right) \in E \times E \mapsto e_j \in E$,
	and consider the exterior differential system
	$p_0^* \om = p_1^* \om$ on $E \times E$.
	This exterior differential system satisfies
	the hypotheses of the Frobenius theorem,
	so integral manifolds have dimension equal
	to the dimension of $E$. Pick any two points $e_0,e_1 \in E$, 
	and let $\Lambda$ be
    the $\Gpt$-orbit of the integral manifold containing 
	$\left(e_0,e_1\right)$; so $\Lambda$
	is an integral manifold, and $\vec{A}$-invariant. 
	Therefore $p_j : \Lambda \to E$ are local diffeomorphisms,
	$\vec{A}$ and $\Gpt$ equivariant, so covering maps.
	So $p_j : \Lambda/\Gpt \to E/\Gpt=M$, a covering map.
	But $M$ is connected and simply connected, and 
	$\Lambda/\Gpt$ is connected by construction, so
	$p_j : \Lambda/\Gpt \to M$ is a diffeomorphism,
	and therefore $p_j : \Lambda \to E$ is a
	diffeomorphism. So $\Lambda$ is the
	graph of an automorphism
	taking $e_0$ to $e_1$. Therefore the automorphism
	group $\Gtot'$ acts transitively on $E$,
	and without fixed points. Pick any $e_0 \in E$,
	identify $g' \in \Gtot' \mapsto g' e_0 \in E$,
	and map $g \in \Gpt \to r_g 1 \in \Gtot'$.
	Check \al{that} this is a homomorphism, because automorphisms
	commute with the right action. \al{Check that this is
        an injection,}
	mapping $\Gpt \to \Gpt'$, and acts simply
	transitively on the fiber of $E$ over $m_0 \in M$,
	so $\Gpt = \Gpt'$.
\end{proof}
\begin{remark}
Warning: not every homogeneous
space with invariant Cartan geometry is a mutation.
To be a mutation, its symmetry group $\Gtot'$
must have the same dimension as that of the model.
A homogeneous space $M=\Gtot'/\Gpt'$ bearing a $\Gtot'$-invariant
Cartan geometry $E \to M$ modelled on $\Gtot/\Gpt$ will frequently
have $\dim G' < \dim G$, with 
curvature varying on $E$, constant
only on the $G'$ orbits inside $E$, and 
need not be complete. For instance,
the hyperbolic plane has a flat projective connection
(i.e. modelled on $\Gtot=\PSL{3,\R{}} \to \Gtot/\Gpt = \RP{2}$)
invariant under hyperbolic isometries, but
is incomplete; it is \emph{not} a mutation of
$\RP{2}$, since the group $G'$ (the projective 
automorphism group of the hyperbolic plane) 
is merely $\PSL{2,\R{}}$, too small.
Similarly, the usual Riemannian metric on $\CP{n}$  
gives a real conformal structure,
which is equivalent to a type of Cartan connection, but not a mutation
of the conformal structure of $S^{2n}$ since the conformal symmetry
group of $\CP{n}$, which is $\mathbb{P}U(n+1)$, is much smaller than $\PO{2n+2,1}$
(the conformal symmetry group of the flat conformal structure on the model $S^{2n}$).
\end{remark}
\begin{corollary}
The automorphism
group of a Cartan geometry modelled
on $\Gtot/\Gpt$ has at most 
$\dim \Gtot$ dimensions, and
has $\dim \Gtot$ dimensions
just when it is a mutation.
\end{corollary}
\subsection{Torsion}
\begin{definition}
If $E \to M$ is a Cartan geometry modelled on $\Gtot/\Gpt$,
and $W$ is a $\Gpt$-representation, let $\Gpt$ act
on $E \times W$ on the right by $(e,w)g=\left(eg,g^{-1}w\right)$,
and let $E \times_{\Gpt} W$ be the quotient by that action;
clearly a vector bundle over $E$.
\end{definition}
\begin{lemma}[Sharpe \cite{Sharpe:1997}]
Define the \emph{1-torsion} $\kappa_1$ to be the image of the
curvature under $\gtot \to \gtot/\gpt$, i.e.
\[
\kappa_1 = \kappa + \gpt \otimes
\Lm{2}{\quot{\gtot}{\gpt}}^* : E \to
\quot{\gtot}{\gpt} \otimes
\Lm{2}{\quot{\gtot}{\gpt}}^*.
\]
The tangent bundle is
\[
TM = E \times_{\Gpt} \left( \quot{\gtot}{\gpt} \right )
\]
and the 1-torsion is a section of \( TM \otimes \Lm{2}{T^*M}.  \)
\end{lemma}
Let $\semiprol{\gpt}{0}=\gpt$ and let $\semiprol{\gpt}{p+1}$ be the
kernel of $0 \to \semiprol{\gpt}{p+1} \to \semiprol{\gpt}{p} \to
\gl{\gtot/\semiprol{\gpt}{p}}$ and
$\prol{\gpt}{p}=\semiprol{\gpt}{p}/\semiprol{\gpt}{p+1}$.  Define
$\kappa_{p+1} = \kappa + \semiprol{\gpt}{p} \otimes
\Lm{2}{\gtot/\gpt}^*: E \to
\gtot/\semiprol{\gpt}{p} \otimes
\Lm{2}{\gtot/\gpt}^*$, so that
$\kappa_1$ is the 1-torsion.  Thus if $\kappa_{q+1}=0$ for all $q <
p$ then $\kappa_{p+1} \in g^{(p-1)} \otimes
\Lm{2}{\gtot/\gpt}^*$.
\begin{lemma}
  If $\kappa_{q+1}=0$ for $q<p$ then $\kappa_{p+1}$ determines a
  section of the bundle $E \times_{\Gpt} \prol{\gpt}{p-1}$, which is a subbundle of $TM
  \otimes \Sym{p}{T^*M} \otimes \Lm{2}{T^*M}$.
\end{lemma}
\begin{proof}
  Sharpe \cite{Sharpe:1997} gives an easy proof of the case $p=0$.
  Induction from there is elementary.
\end{proof}
Similarly we define groups $\semiprol{\Gpt}{0} = \Gpt$ and 
$1 \to \semiprol{\Gpt}{p+1} \to
\semiprol{\Gpt}{p} \to \GL{\gtot/\semiprol{\gpt}{p}}$ and
$\prol{\Gpt}{p}=\semiprol{\Gpt}{p}/\semiprol{\Gpt}{p+1}$.
\subsection{Circles}
\begin{definition}
  Let $\Gtot/\Gpt$ be a homogeneous space.  
  A \emph{circle} is a compact curve in $\Gtot/\Gpt$
  which is acted on transitively by a subgroup of
  $\Gtot$. Without loss of generality,  
  \al{assume} that our circles pass through
  the ``origin'' $\Gpt \in \Gtot/\Gpt$.
  A Lie subgroup $K \subset \Gtot$ is called a \emph{circle subgroup}
  if 
	\begin{enumerate}
	\item its orbit in $\Gtot/\Gpt$ through $\Gpt$ 
  is a circle, and 
	\item $K$ is maximal among subgroups
  with that circle as orbit. 
	\end{enumerate}
	Each circle is a homogeneous 
  space $K/(K \cap \Gpt)$.
  Let $E \to M$ be a Cartan geometry modelled on $\Gtot/\Gpt$. \al{Let $K/(K \cap
  \Gpt)$ be a circle} in $\Gtot/\Gpt$. Let $\mathfrak{k}$ be the Lie 
  algebra of $K$.
  The differential equations $\om = 0 \mod \mathfrak{k}$ are
  holonomic, since the curvature is a semibasic 2-form. Therefore
  these equations foliate $E$, and their projections to $M$ are called
  the \emph{circles} of the Cartan geometry.  If $\gtot$ has a
  $K \cap \Gpt$-invariant complement to $\mathfrak{k}$, 
  then each $K$ circle bears
  a $K/(K \cap \Gpt)$-Cartan connection, determined canonically out of
  $\om$.
\end{definition}
Warning: Sharpe \cite{Sharpe:1997} calls the images in $M$ of the
flows of the $\vec{A}$ vector fields \emph{circles}. For example, he
calls a point of $M$ a circle.  Cartan's use of the term \emph{circle}
in conformal geometry \cite{Cartan:68} and in $CR$ geometry
\cite{Cartan:136,Cartan:136bis} agrees with our use.
For example, in flat conformal geometry on a
sphere, our circles are just ordinary circles, while some of 
Sharpe's circles have one point deleted, and others do not.
The reader may wish to consider the group of projective
transformations of projective space, where one sees clearly 
how the two notions of circle diverge. For more on circles,
see \u{C}ap, Slov{\'a}k and \u{Z}adn{\'i}k \cite{CapSlovakZadnik:2004}.
It will be irrelevant for our purposes as to how one
might attempt to parameterize the circles.
\section{Projective connections}
\subsection{The concept}
  A projective connection is a Cartan geometry modelled
  on projective space $\Proj{n}=\PSL{n+1,\C{}}/G_0$ 
  (with $G_0$ \al{being} the stabilizer of a point of $\Proj{n}$).
  Cartan \cite{Cartan:70} shows that any connection on the tangent bundle
  of a manifold induces a projective connection, with two
  torsion-free connections inducing the same projective connection just
  when they have the same unparameterized geodesics. Not
  all projective connections arise this way\footnote{Those which
  do are often called \emph{projective structures}, although
some authors use the term \emph{projective structure} for
  a projective connection which arises from a torsion-free connection, or
  for a flat projective connection.} but nonetheless
  we can use this as a helpful picture. A Riemannian manifold
  has a projective connection induced by its Levi-Civita
  connection.
  
  The projective connection on (real or complex) projective space is
  complete, because the relevant flows generate the action of the
  projective linear group. Clearly no affine chart is complete,
  because the projective linear group doesn't leave it invariant. 
  A flat torus cannot be complete because it is covered by affine space.  So
  projective completeness is quite strong. No real manifolds with 
  complete real projective connections are known other than flat 
  ones;\footnote{\dots until McKay \cite{McKay:2005} \dots}.
  Corollary~\vref{cor:CompleteGivesKlein} and \al{Theorem}~\vref{thm:ConstantCurv}
  say that the flat complete ones are the sphere and real projective space.
  Kobayashi \& Nagano \cite{KobayashiNagano:1964} \al{conjectures} that complete
  projective connections can only exist on compact manifolds.
  We will prove this, and much more, for holomorphic projective
  connections. Lebrun \& Mason \cite{LeBrunMason:2002} gave a beautiful geometric
  description of the projective structures  
  on $S^2$ for which all geodesics are closed, but have not 
  addressed completeness, which is not obvious even in their examples.\footnote{But
  see McKay \cite{McKay:2005}, where this is proven.}
  \al{Their approach does not apply to general} projective connections,
  even if all of their geodesics are closed.
  It is not clear whether
  completeness of a projective connection implies closed geodesics.\footnote{Once 
  again, McKay \cite{McKay:2005} gives complete examples with geodesics \emph{not}
  closed.}
  Pick a geodesic on the round sphere. We can alter the metric in
  a neighbourhood of that geodesic to look like the metric on
  the torus (picture gluing two spherical caps to a cylinder).
  The resulting real projective structure is incomplete, but
  the limit as we remove the alteration is complete. 
\subsection{Literature}
{\'E}lie Cartan \cite{Cartan:70} introduced the notion of projective
connection; Kobayashi \& Nagano \cite{KobayashiNagano:1964},
Gunning \cite{Gunning:1978} and Borel \cite{Borel:2001} 
provide a contemporary review.
\subsection{Notation}
Let $e_{\mu} \in \Proj{n} (\mu = 0,\dots,n)$ be the $\mu$-th standard
basis vector projectivized.  Let $\Gtot=\PSL{n+1,\C{}}$, $\Gpt$ the
stabilizer of the point $e_0$, $\Ggeod$ the stabilizer of the
projective line through $e_0$ and $e_1$, and $\Gptgeod = \Gpt \cap
\Ggeod$.  To each projective connection on a manifold $M$, Cartan
associates a Cartan geometry (see Sharpe \cite{Sharpe:1997}), which is
a principal right $\Gpt$-bundle $E \to M$ with a Cartan connection
$\om \in \nForms{1}{E} \otimes \gtot$, which we can split up into
\[
\om =
\begin{pmatrix}
  \om^0_0 & \om^0_j \\
  \om^i_0 & \om^i_j
\end{pmatrix},
\]
(using the convention that Roman lower case indices start at 1), with
the relation $\om^0_0 + \om^k_k=0$.  It is
convenient to introduce the expressions
\begin{align*}
  \om^i &= \om^i_0 \\
  \om_i &= \om^0_i \\
  \gamma^i_j &= \om^i_j - \delta^i_j \om^0_0. \\
\end{align*}
\subsection{Structure equations}
There are uniquely determined functions 
\[
K^i_{jk},K_{ijk},K^i_{jkl}
\]
so that the structure equations in table~\vref{tbl:Struc} are satisfied.
\begin{table}
\begin{align*}
  \nabla \om^i &= d \om^i + \gamma^i_j \wedge \om^j \\
  &= \frac{1}{2} K^i_{kl} \om^k \wedge \om^l \\
  \nabla \gamma^i_j &= d \gamma^i_j + \gamma^i_k \wedge \gamma^k_j -
  \left( \om_j \delta^i_k + \om_k \delta^i_j
  \right) \wedge \om^k  \\
  &= \frac{1}{2} K^i_{jkl} \om^k \wedge \om^l \\
  \nabla \om_i &= d \om_i - \gamma^j_i \wedge \om_j \\
  &= \frac{1}{2} K_{ikl} \om^k \wedge \om^l \\
  0 &= K^i_{jk} + K^i_{kj} \\
  0 &= K^i_{jkl} + K^i_{jlk} \\
  0 &= K_{ikl} + K_{ilk}. \\
\end{align*}
\caption{The structure equations of a projective connection}\label{tbl:Struc}
\end{table}
The elements $h \in \prol{\Gpt}{1}$ which fix not only the point $e_0 \in
\Proj{n}$, but also the tangent space at $e_0$, look like
\[
g =
\begin{bmatrix}
  1 & g^0_j \\
  0 & \delta^i_j
\end{bmatrix},
\]
these equations say
\[
r_g^*
\begin{pmatrix}
  0 & K_{jkl} \\
  K^i_{kl} & K^i_{jkl} \\
\end{pmatrix}
=
\begin{pmatrix}
  0 & K_{jkl} \\
  K^i_{kl} & K^i_{jkl} \\
\end{pmatrix}
+ g^0_m
\begin{pmatrix}
  0 & -K^m_{jkl} \\
  0 & K^n_{kl} \left( \delta^m_j \delta^i_n + \delta^i_j \delta^m_n
  \right)
\end{pmatrix}.
\]
\subsection{Flatness}
The canonical projective connection on $\Proj{n}$ is given by $E =
\Gtot \to M = \Proj{n}$, with $\om$ \al{being} the left invariant
Maurer--Cartan 1-form on $\Gtot$; it is flat.
\subsection{Projective connections as second order structures}
Take $E \to M$ the bundle of our projective connection $\om$.  Let
$FM \to M$ be the bundle of coframes of $M$, i.e. the linear
isomorphisms $u : T_m M \to \C{n}$. Following Cartan \cite{Cartan:70},
we can map $E \to FM$ by the following trick: at each point of $E$,
we find that $\left(\om^i\right)$ is semibasic for the map $E \to
M$, so descends to a $\C{n}$-valued 1-form at a point of $M$, giving a
coframing.  A similar, but more involved, process gives a map $E \to
FM^{(1)}$ to the prolongation of $FM$, as follows: Cartan shows that
at each point of $E$, $\gamma^i_j$ is semibasic for the map $E \to
FM$, so descends to a choice of connection at that point of $FM$ (a
connection as a $\GL{n,\C{}}$-bundle $FM \to M$). The bundle $FM^{(1)}
\to FM$ is precisely the bundle of choices of connection, hence a map
$E \to FM^{(1)}$. Finally, Cartan shows that
\[
\xymatrix{
  E \ar[dr] \ar[rr] & &  FM^{(1)} \ar[dl] \\
  &M& }
\]
is a principal right $\Gpt$-bundle morphism.
\subsection{Geodesics}
Cartan points out that the geodesics of a projective connection are
precisely the curves on $M$ which are the images under $E \to M$ of
the flow of the vector field dual to $\om^1$. Denote this vector
field \al{as} $\vec{e}_1$.  Given any immersed curve $\phi : C \to M$, we take
the pullback bundle $\phi^* E \to C$, and inside it we find that the
1-form $\left(\om^i\right)$ has rank 1. Applying the $\Gpt$-action,
we can arrange that $\om^I=0$ for $I>1$.  (Henceforth, capital
Roman letters designate indices $I=2,\dots,n$).  Let $E_C \subset
\phi^* E$ be the subset on which $\om^I=0$, precisely the subset of
$\phi^* E$ on which the map $E \to FM$ reaches only coframes
$\om^i$ with $\om^I$ normal to $C$. Taking exterior derivative
of both sides of the equation $\om^I=0$, we find
$\gamma^I_1=\kappa^I \om^1$ for some functions $\kappa^I : E_C \to
\C{}$.  The tensor $\kappa^I \vec{e}_I \left(\om^1\right)^2$
descends to a section of $\left(\phi^* TM/\phi(TC)\right) \otimes
\left(T^*C\right)^2$, the \emph{geodesic curvature}. A geodesic is
just an immersed curve of vanishing geodesic curvature.
\subsection{Rationality of geodesics}
Henceforth we work in the category of holomorphic geometry
unless explicitly stated otherwise.
\begin{proposition}
  A holomorphic projective connection on a complex manifold is
  complete just when its geodesics are rational curves.
\end{proposition}
\begin{proof}
  Consider the equations $\om^I=\gamma^I_1=0$ on the bundle $E$.
  These equations are holonomic (i.e. satisfy the conditions of the
  Frobenius theorem), so their leaves foliate $E$. We have seen that
  the leaves are precisely the $E_C$ bundles over geodesics $C$. On
  $E_C$ the 1-forms
\[
\om^1,\gamma^1_1,\gamma^1_J,\gamma^I_J,\om_1,\om_J
\]
form a coframing, with the structure equations
\begin{align*}
  d \om^1 &= - \gamma^1_1 \wedge \om^1 \\
  d \gamma^1_1 &= 2 \om_1 \wedge \om^1 \\
  d \gamma^1_J &= - \gamma^1_1 \wedge \gamma^1_J - \gamma^1_K \wedge
  \gamma^K_J
  + \om_J \wedge \om^1 \\
  d \gamma^I_J &= - \gamma^I_K \wedge \gamma^K_J
  + \delta^I_J \om_1 \wedge \om^1\\
  d \om_1 &= \gamma^1_1 \wedge \om_1 \\
  d \om_J &= \gamma^1_J \wedge \om_1 - \gamma^K_J \wedge \om_K \\
\end{align*}
which are the Maurer--Cartan structure equations of $\Ggeod$, hence
$E_C \to C$ is a flat Cartan geometry modelled on $\Ggeod/\Gptgeod$.
Note that $E_C \to C$ is a $\Gptgeod$ bundle, and that
$\Ggeod/\Gptgeod=\Proj{1}$ is simply connected. By
corollary~\vref{cor:CompleteGivesKlein}, if the projective structure
is complete, then $E_C \to C$ is a complete flat Cartan connection, so
$C$ is covered by $\Proj{1}$, an unramified cover, and therefore a
biholomorphism.

Conversely, if the geodesic $C$ is rational, then
theorem~\vref{thm:compactToFlat} proves that $E_C$ is a locally Klein
geometry on $C$, so the vector field $\vec{e}_1$ is complete on $E_C$.
These $E_C$ manifolds foliate $E$, so $\vec{e}_1$ is complete on $E$.
By changing the choice of indices, we can show that
$\vec{e}_2,\dots,\vec{e}_n$ are also complete on $E$.
\end{proof}
\begin{definition}
  A geodesic $C$ is called \emph{complete} when the Cartan geometry on
  $E_C$ is complete.
\end{definition}
\begin{proposition}
  A geodesic is complete just when it is rational; its normal bundle
  is then $\OOp{1}{n-1}$.
\end{proposition}
\begin{proof}
  Let $\phi : C \to M$ be our geodesic.  We have seen that
  completeness is equivalent to rationality.  Split $\C{n}=\C{} \oplus
  \C{n-1}$. Get the structure group $\Gptgeod$ to act on $E_C
  \times \C{n}$ in the obvious manner on $E_C$, and \al{on} $\C{n}$ by
  taking
\[
g=
\begin{bmatrix}
  1 & g^0_1 & g^0_J \\
  0     & g^1_1 & g^1_J \\
  0     & 0     & g^I_J \\
\end{bmatrix},
\]
to act on vectors in $\C{n}$ as the matrix
\[
\begin{pmatrix}
  g^1_1 & g^1_J \\
  0     & g^I_J \\
\end{pmatrix}.
\]
Sharpe \cite{Sharpe:1997} p. 188 gives an elementary demonstration
that $TC = E_C \times_{\Gptgeod} \C{}$, and the same
technique shows that
\begin{align*}
  \phi^* TM &= E_C \times_{\Gptgeod} \C{n}, \\
  \quot{\phi^* TM}{TC} &= E_C \times_{\Gptgeod} \C{n-1}.
\end{align*}
Because $E_C = E_{\Proj{1}}$ are the same Cartan geometries, the
corresponding vector bundles are isomorphic.
\end{proof}
\subsection{Completeness}
\begin{theorem}
  The only holomorphic projective connection on any complex manifold
  which (1) is complete or (2) has all geodesics rational
  is the canonical projective connection on $\Proj{n}$.
\end{theorem}
\begin{proof}
  All geodesics must be rational, with tangent bundles
  $\OO{2}$ and normal bundles $\Norm{C}=\OOp{1}{n-1}$.
  Because $\Gptgeod$ acts trivially on $K^i_{kl}$, the 1-torsion
  (represented by $K^i_{kl}$) lives in the vector bundle $\phi^* TM
  \otimes \Lm{2}{T^*M}$, which is a sum of line bundles of the form
  $\OO{d_i-d_k-d_l}$, where each $d_i,d_k,d_l$ is either 1 or
  2. As long as we take $k=1,$ i.e. only plug in tangent vectors to
  the geodesic to one of the slots in the 1-torsion, we have a section
  of $\phi^* TM \otimes T^*C \otimes \Conorm{C}$.  This maps to a section of
  $\Norm{C} \otimes T^*C \otimes \Conorm{C}$ in which the line bundles are
  negative, so there are no nonzero holomorphic sections. Therefore
  $K^I_{1L}=0$ for $I,L \ne 1$.  This holds at every point of $E_C$,
  and varying the choice of $C$, these $E_C \subset E$ fill out all of
  $E$.  By $\Gpt$-equivariance of $\kappa_1$, we must have
  $\kappa_1=\left(K^i_{kl}\right)=0$ everywhere.
  $\kappa_2=\left(K^i_{jkl}\right)$ becomes a section of the vector
  bundle $\phi^* TM \otimes \Lm{3}{T^*M}$, and $K^I_{11L}=0$ by a
  similar argument, and then $K^I_{J1L}$ becomes defined, and vanishes
  by the same argument, and then finally $\kappa_2=0$ by $\Gpt$
  equivariance. With that established,
  $\kappa_3=\left(K_{jkl}\right)$ becomes a section of $T^*M
  \otimes \Lm{2}{T^*M}$, and in the same way one sees that
  $K_{jkl}=0$.  Apply corollary~\vref{cor:mutual} to see that
  our manifold is $\Gamma \backslash \Proj{n}$ (since $\Proj{n}$
  is simply connected). $\Gamma$ is a discrete group of projective 
  linear transformations, acting freely and properly.  
  Every complex linear transformation has an eigenvalue,
  so every projective linear transformation has a fixed point,
  and therefore $\Gamma$ cannot act freely unless $\Gamma=\{1\}$.
\end{proof}
\section{Conformal connections}
{\'E}lie Cartan \cite{Cartan:70} introduced the notion of a space with
conformal connection (also see Kobayashi \cite{Kobayashi:1995}). \al{These
spaces} include conformal geometries, i.e.  pseudo-Riemannian geometries up to
rescaling by functions.  Since the line of reasoning is so similar
here, it is only necessary to outline the differences between
conformal connections and projective connections.  Again, holomorphic
conformal connections are Cartan connections modelled on $\Gtot/\Gpt$
where $\Gtot$ is the stabilizer of a smooth hyperquadric in projective
space, i.e. the group of complex conformal transformations 
$\Gtot=\PO{n+2,\C{}}$, 
and $\Gpt$ is the stabilizer of a point of the hyperquadric.
Cartan takes the hyperquadric
\[
x_1^2 + \dots + x_n^2 + 2 x_0 x_{n+1} = 0
\]
inside $\Proj{n}$, and \al{takes} $\Gpt$ to be the
subgroup fixing the point $e_0$. Denote the Lie algebras of
$\Gtot$ and $\Gpt$ by the symbols $\gtot$ and
$\gpt$ respectively.  A Cartan connection with this model is an 
$\Gpt$-equivariant differential form $\om \in \nForms{1}{E} \otimes \gtot$
on a principal right $\Gpt$-bundle $E \to M$,
\[
\om =
\begin{pmatrix}
  \om^0_0 & \om^0_j & \om^0_{n+1} \\
  \om^i_0 & \om^i_j & \om^i_{n+1} \\
  \om^{n+1}_0 & \om^{n+1}_j & \om^{n+1}_{n+1}
\end{pmatrix}
\]
subject to
\begin{align*}
  \om^{n+1}_0 &= \om^0_{n+1} = 0 \\
  \om^{n+1}_{n+1} &= - \om^0_0 \\
  \om^i_0 &= -\om^{n+1}_i \\
  \om^i_{n+1} &= - \om^0_i \\
  \om^i_i &= \om^i_j + \om^j_i = 0
\end{align*}
where $i,j=1,\dots,n$. These equations tells us just that $\om$ is
valued in $\gtot$.  It is convenient to define
\begin{align*}
  \varpi   &= \om^0_0 \\
  \om^i &= \om^i_0 \\
  \om_i &= \om^0_i \\
  \gamma^i_j &= \om^i_j - \delta^i_j \om^0_0 \\
\end{align*}
so that
\[
\om = 
\begin{pmatrix}
\varpi & \om_j & 0 \\
\om^i & \gamma^i_j + \delta^i_j \varpi & -\om_i \\
0 & -\om^j & -\varpi
\end{pmatrix}.
\]
\al{Table~\ref{table:confStruc}} gives the structure equations of conformal connections.
\begin{table}
\begin{align*}
  \nabla \varpi &= d \varpi + \om_k \wedge \om^k \\
                &= \frac{1}{2} K_{kl} \om^k \wedge \om^l \\
  \nabla \om^i &= d \om^i + \gamma^i_j \wedge \om^j \\
                &= \frac{1}{2} K^i_{kl} \om^k \wedge \om^l \\
  \nabla \gamma^i_j &= d \gamma^i_j + \gamma^i_k \wedge \gamma^k_j +
  \left( \delta^j_k \om_i - \delta^i_j \om_k -
    \delta^i_k \om_j \right)
  \wedge \om^k \\
  &= \frac{1}{2} K^i_{jkl} \om^k \wedge \om^l \\
  \nabla \om_j &= d \om_j + \om^k \wedge \gamma^k_j \\
                  &= \frac{1}{2} K_{jkl} \om^k \wedge \om^l \\
\end{align*}
\caption{The structure equations of conformal connections}\label{table:confStruc}
\end{table}
Each of the curvature functions $K$ is antisymmetric in the
$k,l$ indices. Elements of $\prol{\Gpt}{1}$ have the form
\[
g =
\begin{bmatrix}
  1 & g^0_j & -\frac{1}{2} g^0_k g^0_k \\
  0 & \delta^i_j & - g^0_i \\
  0 & 0 & 1
\end{bmatrix}
\]
and satisfy the equations in table~\vref{table:pullBackConf}.
\begin{table}
\begin{align*}
  r_g^*
\begin{pmatrix}
  \varpi \\
  \om^i \\
  \gamma^i_j \\
  \om_j \\
\end{pmatrix}
=&
\begin{pmatrix}
  \varpi - g^0_k \om^k \\
  \om^i \\
  \gamma^i_j + \left( \delta^i_k g^0_j - \delta^j_k g^0_i + \delta^i_j
    g^0_k
  \right) \om^k \\
  \om_j - g^0_k \gamma^k_j - \left( g^0_j g^0_k - \frac{1}{2}g^0_l
    g^0_l \delta_{jk} \right) \om^k
\end{pmatrix}
\\
r_g^*
\begin{pmatrix}
  K_{kl} \\
  K^i_{kl} \\
  K^i_{jkl} \\ 
  K_{jkl} \\
\end{pmatrix}
=&
\begin{pmatrix}
  K_{kl} \\ 
  K^i_{kl} \\ 
  K^i_{jkl} \\
  K_{jkl} \\
\end{pmatrix}
\\
&+
\begin{pmatrix}
  -g^0_m K^m_{kl} \\
  0 \\
  \left( \delta^i_m g^0_j - \delta^j_m g^0_i + \delta^i_j g^0_m
  \right) K^m_{kl} \\
  -g^0_m K^m_{jkl} - \left( \delta^i_m g^0_j -
    \delta^j_m g^0_i + \delta^i_j g^0_m
  \right ) K^m_{kl} \\
\end{pmatrix}
\end{align*}
\caption{Action of the structure group on
the curvature of a conformal connection}\label{table:pullBackConf}
\end{table}
\subsection{Circles}
Again, we use capital Roman letters for 
indices $I=2,\dots,n$.  Let
$\Ggeod$ be the subgroup of $\Gtot$ fixing the rational curve
$\Proj{1} = \left(x^I=0\right)$ (this curve is a circle), and
$\Gptgeod=\Gpt\cap\Ggeod$.
\begin{lemma}
  The normal bundle of this curve inside the hyperquadric is
  $\Norm{\Proj{1}} = \OOp{2}{n-1}$.
\end{lemma}
\begin{proof}
  An entirely elementary calculation in projective
  coordinates.
\end{proof}
The group $\Ggeod$ has Lie algebra cut out by the equations
$\om^I=\gamma^I_1=\om_I=0$.  Take $E \to M$ any Cartan geometry
modelled on the hyperquadric $\Gtot/\Gpt$, with Cartan connection
$\om$; Cartan \cite{Cartan:70} shows that the same equations
$\om^I=\gamma^I_1=\om_I$ foliate $E$ by $\Gptgeod$-subbundles
over curves in $M$ which he called \emph{circles}. They are
also circles in our sense.
\begin{theorem}
  The only holomorphic conformal connection on any complex manifold
  which is either complete or has all circles rational
is the canonical flat conformal connection on a 
  smooth hyperquadric.
\end{theorem}
\begin{proof}
  All of the line bundles appearing in the normal and tangent bundles
  of the circles are the same, $\OO{2}$, and therefore 1-torsion lives
  in a sum of $\OO{-2}$ line bundles, and similarly all higher
  torsions live in sums of negative line bundles. Following
  the same line of reasoning as we did for projective connections,
  this implies that our manifold is a quotient of the hyperquadric
  by a discrete subgroup of complex 
  conformal transformations. But every complex conformal
  transformation is orthogonally diagonalisable, so leaves invariant some
  flag of complex subspaces, one of which intersects the
  hyperquadric in a pair of points. Therefore there are no smooth
  quotients of the hyperquadric.
\end{proof}
\section{Rational homogeneous varieties}
\subsection{Locally Klein geometries}
Linear algebra eliminated the possibility of a quotient of
projective space or the hyperquadric by a discrete group $\Gamma$. 
A little Lie theory covers these \al{cases} and many other cases at once.
Recall that a \emph{parabolic} subgroup of a complex semisimple Lie
group is a closed subgroup which contains a Borel (i.e. maximal 
connected solvable) subgroup. See Serre \cite{Serre:2001}
for all of the theory of semisimple Lie algebras that will be used in this article.
\begin{lemma}\label{lemma:GPquotients}
There are no smooth quotients of any complex rational homogeneous variety
$\Gtot/\Parab$ by any discrete subgroup $\Gamma \subset G$ 
with $\Gtot$ \al{being} semisimple and $\Parab$ \al{being} a parabolic
subgroup, i.e. if $\Gamma \backslash G/P$ is a locally
Klein geometry, then $\Gamma=1$.
\end{lemma}
\begin{proof}
Each element $g \in \Gtot$ must belong to a Borel
subgroup (see Borel \cite{Borel:1991} 11.10, \al{p. 151 for proof}),
\al{and all Borel subgroups are conjugate. Let $B$ be a Borel
subgroup.} Every
parabolic subgroup contains a Borel subgroup,
say $h B h^{-1} \subset \Parab$, so $g$ fixes
$h^{-1} P \in \Gtot/\Parab$. Thus no subgroup of $\Gtot$
acts freely on $\Gtot/\Parab$.
\end{proof}
A Cartan geometry modelled on a rational homogeneous
space is called a \emph{parabolic geometry};
\v{C}ap \cite{Cap:2005} explains the theory of 
parabolic geometries beautifully.
\begin{corollary}\label{cor:GP}
Every flat complete holomorphic parabolic geometry 
is isomorphic to its model $\Gtot/\Parab$.
\end{corollary}
\subsection{Homogeneous vector bundles}
\begin{definition}
If $P \subset G$ is a parabolic subgroup of a semisimple Lie group,
and $W$ is a $P$-representation, the vector bundle $G \times_P W \to G/P$
is called a \emph{homogeneous vector bundle}.
\end{definition}
\begin{lemma}
A homogeneous vector bundle $V = G \times_P W$ on a rational
homogeneous variety $G/P$ is equivariantly trivial,
$V = G/P \times W_0$ for $W_0$ a $G$-representation,
just when $W$ is a $G$-representation.
\end{lemma}
\begin{proof}
If $W$ is a $G$-representation, map $(g,w) \in G \times W \to \left(gP,gw\right) \in G/P \times W$.
Check that this is $P$-equivariant, for the $P$-action $(g,w)p=\left(gp,p^{-1}w\right)$
and therefore drops to a map $G \times_P W \to G/P \times W$.
Moreover, this map is obviously equivariant under the left $G$-actions
$g'(g,w)=\left(g'g,w\right)$ on $G \times W$ and $g'\left(gP,w\right)=\left(g'gP,g'w\right)$.

Conversely, given an equivariant isomorphism, $\Phi : G \times_P W \to G/P \times W_0$,
define $\phi : G \to  W^* \otimes W_0$, $\Phi(g,w)P = (gP,\phi(g)w)$.
Using $\phi(1)$ to identify $W$ and $W_0$, the resulting
$\phi$ turns $W$ into a $G$-representation, extending the
action of $P$.
\end{proof}
\subsection{Circles and roots}
Again recall that a \emph{parabolic} subgroup of a complex semisimple Lie
group is a closed subgroup which contains a Borel (i.e. maximal solvable)
subgroup.  All Borel
subgroups are conjugate, so we can stick with a specific choice of one
of them.  The Lie algebra $\parab \subset \gtot$ of a
parabolic subgroup $\Parab \subset \Gtot$ therefore contains a Borel subalgebra
$\mathfrak{b} \subset \gtot$.  With an appropriate choice of
Cartan subalgebra, the Borel subalgebra consists precisely of the sum
of the root spaces of the positive roots.  Every parabolic subalgebra
$\parab$, up to conjugacy, has Lie algebra given by picking a
set of negative simple roots, and letting their root spaces belong to
$\parab$. A compact complex homogeneous space $\Gtot/\Gpt$ with $\Gtot$ \al{being} a
complex semisimple group has $\Gpt=\Parab$ a parabolic subgroup, and
conversely every $\Gtot/\Parab$ is compact, and \emph{is} a rational homogeneous variety;
see Borel \cite{Borel:1991} and Wang \cite{Wang:1954}.  
The negative roots whose root spaces do
not belong to $\parab$ will be called the \emph{omitted roots}
of $\Parab$.
\begin{lemma}
  Suppose that $\Parab \subset \Gtot$ is a parabolic subgroup
  of a semisimple Lie group with Lie algebras $\parab \subset \gtot$. Then each
  omitted root $\alpha$ contains a root vector living in a unique copy
  of $\slLie{2,\C{}} \subset \gtot$.  The associated connected
  Lie subgroup $\SL{2,\C{}}_{\alpha} \subset \Gtot$ projects under $\Gtot \to
  \Gtot/\Parab$ to a rational curve, i.e. a circle. 
\end{lemma}
\begin{proof}
  Picture the root diagram, with the Lie algebra $\mathfrak{b}$ of $B$
  represented by all of the roots on one side of a hyperplane. Take
  any omitted root $\alpha$. It lies on the negative side of the
  hyperplane.  The line through that root gives a subalgebra
  isomorphic to $\slLie{2,\C{}}$. We can see that this projects to a
  curve $C_{\alpha} \subset \Gtot/\Parab$ because the root space of $\alpha$ is
  semibasic for the submersion. This curve is acted on (faithfully, up
  to a finite subgroup) by the subgroup coming from that
  $\slLie{2,\C{}}$ Lie subalgebra.  By Lie's classification of group
  actions on curves (see Sharpe \cite{Sharpe:1997}), 
  the only connected Lie group with Lie algebra
  $\slLie{2,\C{}}$ which acts on a curve is $\PSL{2,\C{}}$ (up to
  taking a covering group), acting on a rational curve. So we will
  write $C_{\alpha}$ as $\Proj{1}_{\alpha}$. The various roots
  $\alpha$ that we can pick to play this role produce a basis for the
  tangent space of $\Gtot/\Parab$.
  
  There is a very slight subtlety: the subgroup we have called
  $\SL{2,\C{}}_{\alpha}$ might actually be isomorphic to $\PSL{2,\C{}}$,
  so we replace it by its universal cover (largely to simplify
  notation), which will then be an immersed subgroup $\SL{2,\C{}}_{\alpha} \to G$,
  not necessarily embedded.
\end{proof}
  If $C$ is any circle, it must be the orbit of a subgroup $K \subset G$.
  By Lie's classification of group actions on curves,
  the solvable part of $K$ cannot act transitively on $C$, while some 
  subgroup of the semisimple part of $K$ must act transitively and
  must be isomorphic to $\SL{2,\C{}}$ or $\PSL{2,\C{}}$.  
\subsection{Homogeneous vector bundles on the projective line}
Let $B \subset \Sl{}$ be the Borel subgroup of upper triangular matrices, and
$N \subset B$ the maximal nilpotent subgroup of strictly upper triangular matrices.
A homogeneous vector bundle $V$ is then obtained from each
representation $W$ of $B$ by $V = \Sl{} \times_B W$. Moreover, it
is $\Sl{}$-equivariantly trivial just when $W$ is actually a
$\Sl{}$-representation. Let us examine the representations of
$B$. Each one induces a representation of $A \subset B$,  
the maximal reductive subgroup of diagonal matrices. In
particular, the bundle $\OO{-1} \to \Proj{1}$ is given
by taking the representation 
\[
\begin{pmatrix}
a & b \\
0 & a^{-1} 
\end{pmatrix}
w = aw,
\] 
$w \in \C{}$. More generally, the line bundle $\OO{d} \to \Proj{1}$
is given by the representation 
\[
\begin{pmatrix}
a & b \\
0 & a^{-1}
\end{pmatrix}
w = a^{-d} w.
\]
We call a representation $W$ of $B$ \emph{elementary} if
the induced $A$-representation is completely
reducible. An elementary representation on which $N$
acts trivially is precisely a sum of 1-dimensional representations,
and the vector bundle $V$ is a sum of line bundles.
More generally, any elementary representation
$W$ will give rise to a Lie algebra representation
of $\mathfrak{b} = \mathfrak{a} \oplus \mathfrak{n}$ 
on $W$, and this is determined by the matrix of integers which the element
\[
H = 
\begin{pmatrix}
1 & 0 \\
0 & -1
\end{pmatrix} \in \mathfrak{a}
\]
gets mapped to, say 
\[
\begin{pmatrix}
k_1 \\
& k_2 \\
& & \ddots \\
& & & k_n
\end{pmatrix}
\]
and the matrix $\rho(X)$ which the element
\[
X = 
\begin{pmatrix}
0 & 1 \\
0 & 0 
\end{pmatrix}
\in \mathfrak{n}
\]
gets mapped to, say 
\[
\rho(X) = \begin{pmatrix} x^i_j \end{pmatrix}.
\]
Checking the commutation relations of $X$ and $H$, we find
that for each $i,j$, either $k_i-k_j = 2$ or $x^i_j=0$. Therefore 
we can divide up
the indices $i=1,\dots,n$ into partitions, each partition
\al{being a maximal sequence of indices,} so that (1) each index $i$ 
in a sequence has associated 
value $k_i$ larger by 2 than the last index, and (2) if $i$ 
and $j$ are indices which are adjacent in a sequence, then 
$x^i_j \ne 0$. The partitions will be called \emph{strings}
(loosely following Serre \cite{Serre:2001}). Thus every representation
$W$ is a sum of indecomposable representations $W = \sum_s W_s$,
a sum over strings. We draw the strings along a number line,
placing vertices at the values $k_i$. As matrices, once
we organize into a sum of strings,
\[
\rho(H)
=
\begin{pmatrix}
k_1 \\
 & k_1 - 2 \\
 & & \ddots \\
 & & & k_2 \\
 & & & & k_2 -2 \\
 & & & & & \ddots  
\end{pmatrix},
\ \rho(X)
=
\begin{pmatrix}
0  & x_1 \\
  & 0 & x_2 \\
  & & \ddots & \ddots
\end{pmatrix}.
\]

Change to an equivalent representation, say $\rho'=t \rho t^{-1}$
for some invertible matrix $t$, and you obtain an isomorphic
vector bundle. We want to keep $\rho'(H)$ diagonalized in strings, 
so we can permute the blocks coming from the various strings,
and scale the eigendirections, say by
\[
t =
\begin{pmatrix}
\lambda_1 \\
& \lambda_2 \\
& & \ddots
\end{pmatrix}
\] 
This rescales these $x_i$ by $x_i' = \lambda_i x_i \lambda_{i+1}^{-1}$.
Therefore we can arrange the $x_i$ of each string to be anything other than $0$.
In particular, for a single string that has 
\[
\rho(H)
=
\begin{pmatrix}
k \\
& k-2 \\
& & \ddots
\end{pmatrix}
\]
we will always choose to scale to arrange $x_1=k,x_2=k-1,\dots,x_j=k+1-j,\dots$,
\[
\rho(X) = 
\begin{pmatrix}
0 & k \\
 & 0 & k-1 \\
 & & \ddots & \ddots
\end{pmatrix}.
\]
\begin{lemma}
Every indecomposible elementary representation $W$ of $B$, up to 
isomorphism, is obtained by taking a sequence of integers
$k,k-2,\dots,k-2n$, \al{for some} $n \ge 0$, arranging them as the
diagonal entries of a matrix $\rho(H)$,
and arranging the numbers $k,k-1,\dots,$ above the diagonal
of a matrix $\rho(X)$ whose other entries all vanish.
This representation arises from a representation
of $\Sl{}$ just when $k=n$, i.e. when the matrix $\rho(H)$
has entries $k,k-2,\dots,-k$.
In particular, the vector bundle $V = \Sl{} \times_B W$
is $\Sl{}$-equivariantly trivial just when $k=n$.
\end{lemma}
\begin{proof}
Clear from Serre \cite{Serre:2001} p. 19,
taking the standard representations of $\SL{2,\C{}}$.
\end{proof}
\begin{example}
The vector bundle $\OO{-1}$ is drawn in \al{Figure~\ref{fig:OMinusOne}},
\begin{figure}
\[
\xymatrix{
\dots \ar@{-}[rrrrrrrr] & {}_| \ar@{}[]_{-2} & {}_| \ar@{}[]_{-1} & {}_| \ar@{}[]_{0} & \bullet \ar@{}[]_{1} & 
{}_| \ar@{}[]_{2} &
{}_| \ar@{}[]_{3}
& {}_| \ar@{}[]_{4} &
\dots \\
}
\]
\caption{}\label{fig:OMinusOne}
\end{figure}
showing that the string is really just a single point. On the other hand,
the representation $W=\slLie{2,\C{}}$ gives a trivial vector bundle,
and in a picture, it gives the string
in \al{Figure~\ref{fig:SLTwoBundle}}
\begin{figure}
\[
\xymatrix{
\dots \ar@{-}[rrrrrrrr] 
& \bullet \ar@{}[]_{-2} \ar@/^/[rr] 
& {}_| \ar@{}[]_{-1} 
& \bullet \ar@{}[]_{0} \ar@/^/[rr] 
& {}_| \ar@{}[]_{1} 
& \bullet \ar@{}[]_{2} 
& {}_| \ar@{}[]_{3}
& {}_| \ar@{}[]_{4}
& \dots \\
}
\]
\caption{}\label{fig:SLTwoBundle}
\end{figure}
Therefore the tensor product of these vector bundles is given by
the string in \al{Figure~\ref{fig:Tensor}}.
\begin{figure}
\[
\xymatrix{
\dots \ar@{-}[rrrrrrrr] 
& {}_| \ar@{}[]_{-2}
& \bullet \ar@{}[]_{-1}  \ar@/^/[rr] 
& {}_| \ar@{}[]_{0} 
& \bullet \ar@{}[]_{1}  \ar@/^/[rr] 
& {}_| \ar@{}[]_{2} 
& \bullet \ar@{}[]_{3} 
& {}_| \ar@{}[]_{4}
& \dots \\
}
\]
\caption{}\label{fig:Tensor}
\end{figure}
\end{example}
Every $B$-invariant subspace of a string is obtained by removing
some nodes on the left side of the string. The quotient by that 
subspace is obtained by removing the complementary nodes from the right side.
\begin{lemma}
A $B$-representation $W$ is elementary just when it is 
isomorphic to a quotient $W=U/U_0$ of an $\SL{2,\C{}}$-representation
$U$ by a $B$-invariant subspace $U_0$.  Clearly every
elementary representation is a sum of strings.
\end{lemma}
\begin{example}
The vector bundle $\OO{-1} \otimes \left(\SL{2,\C{}} \otimes_B \left(\slLie{2,\C{}}/\mathfrak{b}\right)\right)$ is given by
the string in \al{Figure~\ref{fig:OOne}}
\begin{figure}
\[
\xymatrix{
\dots \ar@{-}[rrrrrrrr] 
& {}_| \ar@{}[]_{-2}
& \bullet \ar@{}[]_{-1} 
& {}_| \ar@{}[]_{0} 
& {}_| \ar@{}[]_{1} 
& {}_| \ar@{}[]_{2} 
& {}_| \ar@{}[]_{3} 
& {}_| \ar@{}[]_{4}
& \dots \\
}
\]
\caption{}\label{fig:OOne}
\end{figure}
and is therefore $\OO{1}$.
\end{example}
\begin{lemma}
The homogeneous vector bundles obtained from elementary representations are precisely those
which are sums of representations of the form $\OO{d} \otimes V$, where $V$ is equivariantly trivial.
\end{lemma}
\begin{proof}
After tensoring with $\OO{-d}$ for an appropriate choice of $d$, we can arrange
that the string sits symmetrically about the origin,
so the same as for an $\SL{2,\C{}}$-representation. Keep in mind that
in the process we have to rescale the various $x_i$ discussed above,
to ensure that they match the $\s{}$-representation.
\end{proof}
\subsection{Normal bundles of circles}
\begin{lemma}
Let $P \subset G$ be a parabolic subgroup of a 
semisimple Lie group $G$ with Lie algebras
$\parab \subset \gtot$.
\al{Every circle $C \subset G/P$, modulo $G$ action, is given as
follows. A root $\alpha$ omitted from $P$ is chosen.
The subalgebra $\s{\alpha} \subset \gtot$
given as
\[
\s{\alpha} = \gtot^{\alpha} \oplus \gtot{-\alpha} \oplus \left[\gtot^{\alpha},\gtot^{-\alpha}\right],
\]
is constructed. The 
connected Lie subgroup of $G$ with this Lie algebra
is constructed. Its universal covering group, say $\Sl{\alpha} \to \Gtot$
is constructed. The subalgebra $\mathfrak{b}_{\alpha} \subset \s{\alpha}$,
\[
\mathfrak{b}_{\alpha} = \gtot^{-\alpha} \oplus \left[\gtot^{\alpha},\gtot^{-\alpha}\right],
\]
is formed. Denote the associated connected subgroup as 
$B_{\alpha} \subset \Sl{\alpha}$. Finally let $C = \Proj{1}_{\alpha} = \Sl{\alpha}/B_{\alpha}$.}
\end{lemma}
\begin{definition}
If $\alpha$ is a root of a Lie algebra $\gtot$, an $\alpha$-\emph{string} 
is a maximal sequence of roots $\beta,\beta-\alpha,\dots$ of $\gtot$
(including the possibility that any one of these roots is allowed to be $0$).
\al{The \emph{length} of a string is the number
of its nodes (unrelated to the Killing metric).}
\end{definition}
\begin{lemma}
Let $P \subset G$ be a parabolic subgroup of a semisimple Lie group, with Lie algebras
$\parab \subset \gtot$.
If $\Proj{1}_{\alpha}=\Sl{\alpha}/B_{\alpha} \subset G/P$ is a circle,
then 
\[
\left.T\left(\quot{G}{P}\right) \right|_{\Proj{1}_{\alpha}} =
 \bigoplus_s \OOp{d_s}{n_s},
\]
a sum over strings, where the strings are precisely
the $\alpha$-strings for $\gtot$, say 
\[
s = \beta,\beta-\alpha,\beta-2 \alpha,\dots,
\]
the rank $n_s$ is the number of omitted roots
in the string $s$, and the degree $d_s$ is the number
of roots in the string $s$ which are not omitted.
\end{lemma}
\begin{proof}
The tangent bundle is $T(G/P)= G \times_P \left(\gtot/\parab\right)$.
Each $\alpha$-string is clearly a representation of $\Sl{\alpha}$,
so gives a trivial vector bundle. The roots that arise from
$\gtot/\parab$ are just the omitted ones. Shift them over by line bundles
to make them into trivial bundles, i.e. shift over the
strings so that they are symmetric under reflection by
$\alpha$. Clearly the number of steps needed is $d_s$, since with this
piece added the result is invariant under $\alpha$ reflection.
\end{proof}
\begin{definition}
The \emph{curvature bundle} on a rational homogeneous
variety $G/P$ is the bundle
\[
G \times_P 
\left(
      \gtot \otimes \Lm{2}{\quot{\gtot}{\parab}}^*
      \right)
\]
Similarly, the curvature bundle on a parabolic geometry
$E \to M$ modelled on $G/P$ is
\[
E \times_P \left(
      \gtot \otimes \Lm{2}{\quot{\gtot}{\parab}}^*
      \right).
\]
\end{definition}
\begin{corollary}
If $\Proj{1}_{\alpha}=\Sl{\alpha}/B_{\alpha} \subset G/P$ is a circle,
then the curvature bundle restricted to $\Proj{1}_{\alpha}$, say
\[
\Sl{\alpha} \times_{B_{\alpha}}
\left(
      \gtot \otimes \Lm{2}{\quot{\gtot}{\parab}}^*
      \right)
\]
is isomorphic to 
\[
\al{
\bigoplus_{s} 
\OOp{0}{\left(n_s+d_s\right)} 
\otimes 
\Lm{2}%
{
	\bigoplus_{s'} 
	\OOp{-d_{s'}}{n_{s'}}
},
}
\]
a sum over $\alpha$-strings $s,s'$ in $\gtot$.
The sections of the curvature bundle
defined on $\Proj{1}_{\alpha}$
span a subbundle, isomorphic to the same sum,
but with $s'$ restricted to $\alpha$-strings
containing no $\parab$-roots.
\end{corollary}
\begin{corollary}
The  sections $K$ of the curvature
bundle which are defined everywhere on $\Proj{1}_{\alpha}$
satisfy $v \hook K = 0$ for $v \in T \Proj{1}_{\alpha}$.
\end{corollary}
\begin{proof}
Mapping 
\[
K \in \gtot \otimes 
\Lm{2}{\quot{\gtot}{\parab}}^* 
\mapsto \bar{K} \in \gtot \otimes 
(\quot{\gtot}{\parab})^* \otimes 
\left(
	\quot{\s{\alpha}}%
	{\mathfrak{b}_{\alpha}}
\right)^*,
\]
we are stripping off all of the possible $\alpha$-strings from the third factor
except for the $\alpha$-string $s'=\alpha,0,-\alpha$, which has $d_{s'}=2$
since $0,-\alpha$ are nonnegative roots so belong to $\parab$. Therefore
$\bar{K}$ lies in a sum of negative degree line bundles, so vanishes.
\end{proof}
\begin{corollary}
The curvature bundle on a rational homogeneous variety
has no global sections.
\end{corollary}
\begin{proof}
The circles issue forth in a basis of directions.
\end{proof}
\section{Complete parabolic geometries}
\begin{theorem}
A parabolic geometry is complete just when all of its
circles are rational.
\end{theorem}
\begin{proof}
Let $E \to M$ be a parabolic geometry.  
If complete, then all of the circles can be developed
onto those of the model, as we saw for projective and
conformal connections, and therefore must be rational.
Conversely, if all circles are rational, then the developing
onto the model circles gives rise to no problems
of monodromy, and clearly the circles 
have induced flat Cartan connections modelled on the circles
of the model. The circles are rational,
hence simply connected, and are flat. By \al{Theorem}~\vref{thm:compactToFlat} they are
complete. The flows of the various
$\vec{A}$ vector fields on $E$ for $A \in \s{\alpha}$
are tangent to the submanifolds $E|_C$ for $C$ \al{being} any
$\alpha$-circle, and on each $E|_C$ they arise
from the induced Cartan geometry. Therefore the
$\vec{A}$ on $E$ are complete.
\end{proof}
\begin{theorem}
In a complete parabolic geometry $E \to M$ modelled
on a rational homogeneous variety $G/P$, every circle $C \subset M$ 
modelled on some $\Proj{1}_{\alpha} \subset G/P$ has
ambient tangent bundle 
$\left.TM\right|_C = \left.T\left(G/P\right)\right|_{\Proj{1}_{\alpha}}$
identified by the developing map.
\end{theorem}
\begin{proof}
The induced Cartan geometries $E|_C \to C$ and 
$G|_{\Proj{1}_{\alpha}} \to \Proj{1}_{\alpha}$ 
are isomorphic by developing (as for projective and conformal connections).
The ambient tangent bundle is $TM|_C = E|_C \times_P \left(\gtot/\parab\right)$,
determined by the induced Cartan geometry.
\end{proof}
\begin{theorem}
Complete parabolic geometries are isomorphic to their models,
i.e. the only complete parabolic geometry modelled on $G/P$
is the standard flat one on $G/P$, modulo isomorphism.
\end{theorem}
\begin{proof}
Let $E \to M$ be our parabolic geometry, modelled
on $G \to G/P$, with $P \subset G$ having Lie algebras $\parab \subset G$.
Let $C$ be a circle modelled on $\Proj{1}_{\alpha}$.
If complete, then all circles are rational. Thus each circle
$C$ has ambient tangent bundle a sum
\[
TM|_C = \bigoplus_s \OOp{d_s}{n_s}
\]
over $\alpha$-strings in the root system of $\gtot$.
The curvature will then live in the \al{curvature bundle:}
\[
\al{
\bigoplus_{s} 
	\OOp{0}{n_s+d_s} 
	\otimes 
	\Lm{2}
	{
		\bigoplus_{s'} 
		\OOp%
		{%
			-d_{s'}%
		}%
		{%
			n_{s'}%
		}%
	},
}
\]
summing over $\alpha$-strings $s$ and $s'$.
If all of the $d_{s'}$ are positive, then
there are no global sections of this
bundle. More generally, global sections always
live in the subbundle where $d_{s'}=0$, a trivial
bundle. One sees immediately that this subbundle
is invariantly defined, since it is spanned
by the global sections of the curvature bundle.  

Any global
section of the curvature bundle on $E \to M$
will be defined on every $\alpha$-circle,
and therefore by isomorphism of the 
induced Cartan geometry on that $\alpha$-circle
with the geometry on the model circle $\Proj{1}_{\alpha} \subset G/P$,
we induce a global section of the curvature
bundle restricted to $\Proj{1}_{\alpha}$ (perhaps
not defined on all of $G/P$). But this
already forces the curvature to vanish in the
tangent direction of the $\alpha$-circle.
Therefore curvature on $E$ vanishes in all directions.
\end{proof}
\begin{theorem}\label{thm:Immutable}
Rational homogeneous varieties are immutable.
\end{theorem}
\begin{proof}
Any mutation is complete, and therefore isomorphic
to the model.
\end{proof}
\begin{corollary}
Real rational homogeneous varieties are immutable.
\end{corollary}
\begin{proof}
Complexify the mutation.
\end{proof}
\section{Flag variety connections}\label{sec:FVC}
To see the normal bundle of a rational circle
in a concrete example, consider a complete Cartan connection $\om$ on $E \to M$, modelled
on a partial flag variety, i.e.  $\Gtot=\SL{n,\C{}}$ and $P=\Gpt$ 
is the
stabilizer of a partial flag in 
$0 \subset V_{k_1} \subset \dots \subset V_{k_1 + \dots + k_p} \subset \al{=} \C{n}$.  
Write
 \[
 \om =
 \begin{pmatrix}
 \om^i_j
 \end{pmatrix}.
 \]
 The subalgebra $\parab$ of $P$ sits inside $\gtot$ as
 diagonal blocks of size $k_1 \times k_1, k_2 \times k_2, \dots$, together
 with all elements above those blocks.
 For any fixed choice of $i>j$ for which $\omega^i_j$ is
 not in such a block, look at the Lie subalgebra of
 $\slLie{n,\C{}}$ for which all entries vanish except for those with
 either the upper or lower index, or both, equal to $i$ or to $j$, i.e.
 the 4 corners of a square inside the matrix, symmetric under
 transpose.  The associated group has a rational orbit in the partial flag
 variety.  The orbit is, up to conjugation, just given by fixing all
 subspaces in a complete flag, except for one. These rational circles stick out
 in ``all directions'' at each point of the partial flag variety, i.e. their
 tangent lines span the tangent space of the partial flag variety. (Perhaps
 surprisingly, the parabolic subgroup does not have any open orbits when it
 acts on tangent lines to the partial flag variety at a chosen point.)
 
\begin{figure}
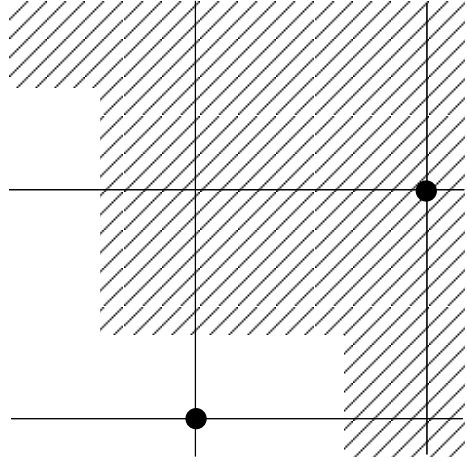

   \centering 
\drawfile{flagvariety}
 \caption{$A_n=\slLie{n+1,\C{}}$: Drawing the 
   matrix $\om$, the unshaded parts of the vertical and horizontal
   lines through the unshaded (lower left) dot in the matrix contribute
   degree 1 line bundles to the normal bundle of the circle; 
   the unshaded dot contributes degree 2 to the tangent
   bundle, and nothing to the normal bundle.}\label{fig:grassmannian}
 \end{figure}

 Developing these rational circles into the Cartan geometry gives a
 family of rational circles in $M$. Consider the rational circle in the
 $\om^i_j$ direction.  Think of the matrix
 $\om$ as having the parabolic subalgebra on and above the diagonal,
 and in the blocks down the diagonal.
 Then the entries
 below the diagonal can be thought of as the tangent space of
 $\Gtot/B$. The $\om^i_j$ entry is the direction of our rational
 curve, which has tangent bundle $\OO{2}$. 
 We think of each choice of two elements $i,j$ as determining
 a choice of root of $G$, the root $\omega^j_j - \omega^i_i$.
 The associated root space is spanned by the matrix whose entries
 all vanish except for $\omega^i_j$, i.e. the root spaces are the
 entries of the matrix $\omega$ off the diagonal. Given
 a root $\alpha=\omega^j_j - \omega^i_i$, and a root $\beta=\omega^{j'}_{j'} - \omega^{i'}_{i'}$,
 their difference $\beta-\alpha$ is also a root just when $i=i'$ or $j=j'$, and
 the difference $\beta-2\alpha$ is also a root just when $\alpha=\beta$.
 Therefore the entries in the same row
 and column, below the diagonal blocks, contribute $\OO{1}$ line bundles to
 $T(G/P)_{\Proj{1}_{\alpha}}$, while $\omega^i_j$ contributes a single $\OO{2}$
 line bundle; see figure~\vref{fig:grassmannian}.
 The ambient tangent bundle on the rational circle coming from $\om^i_j$ is
 \[
 T\left(G/P\right)_{\Proj{1}_{\alpha}} = \OO{2} \oplus \OOp{0}{n_0} \oplus
 \OOp{1}{n_1},
 \]
 where $\alpha = \omega^j_j - \omega^i_i$ and $s$ and $t$ are chosen so that
\begin{align*}
k_1 + \dots \al{+} k_{s-1} &< j \le k_1 + \dots \al{+} k_s, \\
k_1 + \dots \al{+} k_{t-1} &< i \le k_1 + \dots \al{+} k_t, \\
\end{align*}
and it then follows that the ranks are
\begin{align*}
n_1 &= n+k_{s+1}+\dots+k_{t-1}-1, \\
n_0 &= \frac{1}{2} \left( n^2 - \left(k_1^2 + \dots \al{+} k_p^2 \right) \right) - n_1.
\end{align*}
\section{Spinor variety connections}
The examples of normal bundles of rational circles might be important
for numerous applications, so we give another example.
Consider the two spinor varieties (see Fulton \& Harris
\cite{FultonHarris:1991} for their definition). It is enough to check
only one of them, say $\SpinorMinus{\C{2n}}$, since the geometry of
circles and their normal bundles is invariant under the outer
automorphism.  The pattern here is similar again.
\[
\om =
\begin{pmatrix}
 \gamma & \xi \\
 \eta & -\trans{\gamma}
\end{pmatrix},
\]
with
\[
0 = \eta + \trans{\eta} = \xi + \trans{\xi},
\]
and the subalgebra $\parab$ given by $\eta=0$ generates a
subgroup $\Parab$, so that $\Gtot/\Parab=\SO{2n,\C{}}/\Parab$ 
is the spinor variety
$\SpinorMinus{\C{2n}}$.  Transposition of matrices intertwines the theory of
$\SpinorMinus{\C{2n}}$ with that of $\SpinorPlus{\C{2n}}$.
 
\begin{figure}
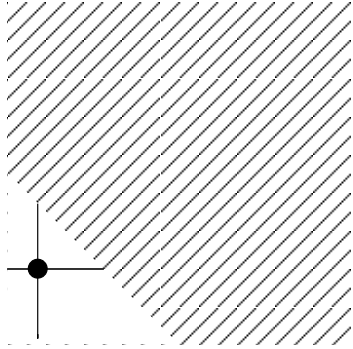

 \centering\drawfile{spinorVariety}
\caption{The spinor variety: drawing the 
 matrix $\om$, the vertical and horizontal lines through the lower
 left dot in the matrix contribute degree 1 line bundles to the
 normal bundle. The other unshaded stuff contributes degree 0.}%
\label{fig:spinorvariety}
\end{figure}

A circle subalgebra is given by setting all but one of the components
of $\eta$ to zero, say leaving only $\eta^i_j$, for some $i>j$,
setting all but two of the components of $\gamma$ to zero, leaving
only $\gamma^i_i$ and $\gamma^j_j$, and setting all but one of the
components of $\xi$ to zero, leaving only $\xi^i_j$, and finally
setting $\gamma^j_j=\gamma^i_i$. The circle subalgebra is isomorphic
to $\slLie{2,\C{}}$.  Its orbit in the spinor variety is therefore a
rational circle.  
The roots of $\SO{2n,\C{}}$ are $\pm \gamma^i_i \pm \gamma^j_j$
for $i \ne j$. The roots $\gamma^i_i + \gamma^j_j$ are precisely
those whose root spaces live in $\gamma=\xi=0$, i.e. the
directions of the tangent space $\gtot/\parab$. Given any
one of these roots, say $\alpha = \gamma^i_i + \gamma^j_j$,
with $i > j$, 
and any other such root $\beta = \gamma^{i'}_{i'} + \gamma^{j'}_{j'}$,
with $i' > j'$, the associated string is
\[
\beta,\beta-\alpha,\beta-2 \alpha, \dots = \gamma^{i'}_{i'} + \gamma^{j'}_{j'}, \gamma^{i'}_{i'} + \gamma^{j'}_{j'} - \gamma^i_i - \gamma^j_j, \dots
\]
and clearly we run out of roots immediately unless $i=i'$ or $i=j'$ or $j=i'$ or $j=j'$.
In such cases, we only obtain a root belonging to $\parab$. 
Moreover, we cannot make another step along the string staying
on a root unless $i=i'$ and $j=j'$. We read off the ambient tangent
bundle of the rational circle in the spinor variety:
\[
T(G/P)_{\Proj{1}_{\alpha}}=\OOp{0}{p_0} \oplus \OOp{1}{p_1}
\]
where
\begin{align*}
 p_0 &= \frac{n^2-3n+6}{2}+j-i \\
 p_1 &= n+i-j-3
\end{align*}
(see \al{Figure}~\vref{fig:spinorvariety}).  
\section{Lagrangian--Grassmannian connections}\label{sec:LagGrass}
For $\Gtot=\Symp{2n,\C{}}$
the Maurer--Cartan 1-form is
\[
\om = 
\begin{pmatrix}
\gamma & \xi \\
\eta & -\trans{\gamma}
\end{pmatrix}
\]
where $\trans{\eta}=\eta$ and $\trans{\xi}=\xi$.
The subalgebra $\parab$ given
by $\om=0$ is the Lie subalgebra of 
the subgroup $\Parab$ preserving a Lagrangian
$n$-plane, so $\Gtot/\Parab$ is the 
\al{Lagrangian--Grassmannian},\mn{Keep this as ``Lagrangian--Grassmannian'',
not ``Lagrangian--Grassmannian connection''.} $\Lag{\C{2n}}$, i.e. the
space of Lagrangian $n$-planes in $\C{2n}$.
Thus $\dim \Lag{\C{2n}}=\frac{n(n+1)}{2}.$
 
Consider the circle subgroup given by
setting all of the $\eta$ and $\xi$
components to $0$, except for one $\eta^i_j$
(which equals $\eta^j_i$) for some
choice of $i \ge j$, the corresponding
$\xi^i_j$, and $\gamma^i_i$ and $\gamma^j_j$,
and finally impose the equation
$\gamma^j_j=\gamma^i_i$. The circle
subalgebra is isomorphic to $\slLie{2,\C{}}$,
and $\eta^i_j$ is semibasic,
\al{and therefore the orbit of this circle
subgroup is 
a rational circle on the Lagrangian--Grassmannian.}
 
The roots of $\Symp{2n,\C{}}$ are $\pm \gamma^i_i \pm \gamma^j_j$
for $i \ne j$ and $\pm 2 \gamma^i_i$. The roots omitted from
$\parab$ are $\gamma^i_i + \gamma^j_j$ (where $i$ might equal $j$).

\begin{figure}
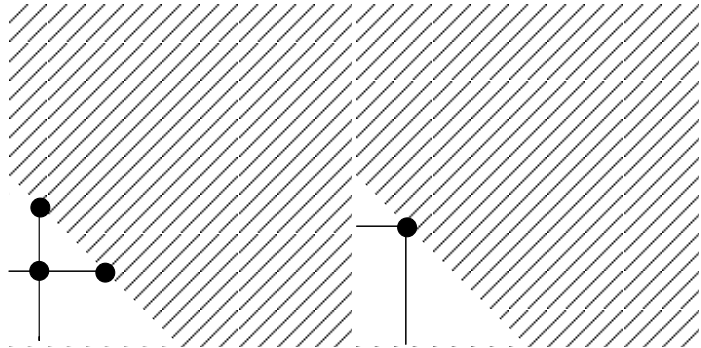

 \centering\subfigure{\drawfile{lag}}%
 \subfigure{\drawfile{lag2}}
\caption{The Lagrangian--Grassmannian: 
 \textbf{(a)} the 
 vertical and horizontal lines
 through the lower left dot in the matrix contribute
 degree 1 line bundles to the normal bundle. 
\textbf{(b)} The count is slightly different 
 when the root space lies on the diagonal.}
\label{fig:lag}
\end{figure}

 The
 components of $\eta$ on or under the diagonal
 (see Figure~\vref{fig:lag}),
 on the same row or column as $\eta^i_j$,
 contribute $\OO{1}$ to the normal
 bundle, while all other $\eta$
 components contribute $\OO{0}$, so
 $T(G/P)_{\Proj{1}_{\alpha}}=\OO{2} \oplus \OOp{0}{p_0} \oplus \OOp{1}{p_1}$
 with
 \[
 p_0 = \frac{n^2+n}{2}-p_1, \quad p_1 = n+i-j-1.
 \]
\section{Conclusions}
\al{Kobayashi invented invariant metrics for
projective structures.} For a generic Cartan geometry on a complex
manifold, the circles of the model roll onto the manifold, but
this rolling runs into singularities, 
rolling out not a rational curve but only a disk.
The hyperbolic metrics on these disks are responsible for
the Kobayashi metric on the manifold: define the distance
between points to be the infimum distance measured along
paths lying on finitely many of these disks, in their 
hyperbolic metrics.

On the other hand, the Cartan geometries attended to in
the present article have quite degenerate Kobayashi
metrics, at the opposite extreme. The next natural step is to
investigate parabolic geometries which are just
barely incomplete, for instance with all circles
affine lines or elliptic curves. Such geometries might
admit classification, although their analysis
seems difficult in the absence of invariant metrics.

The approach of Jahnke \& Radloff \cite{Jahnke/Radloff:2004}
should uncover the sporadic nature of curved parabolic geometries on smooth complex 
projective varieties. Roughly speaking, parabolic geometries appear
to be rare, except perhaps on unspeakably complicated complex manifolds,
and all known parabolic geometries on smooth projective 
varieties are flat. However, there
are simple curved examples with singularities 
(see Hitchin \cite{Hitchin:1982} \al{p. 83, example 3.3}), which
would appear to point the only way forward for
research in parabolic geometry: curvature
demands singularity.\mn{Bibliography entries updated where possible.}
\bibliographystyle{amsplain}
\bibliography{completeness}
\end{document}